\documentclass[10pt,a4paper,oneside,headinclude,footinclude,BCOR5mm]{article}

\usepackage{fancyhdr,amssymb,amsbsy,amsthm}

\title{\normalfont{The $\beta$-symplectic critical surfaces in Hermite surfaces}}

\author{Yongpin Zhu\textsuperscript{1}}

\date{}

\begin{document}


\renewcommand{\sectionmark}[1]{\markright{#1}} 

\pagestyle{fancy} 


\maketitle 






\section*{Abstract}

In \cite{Zhu}, the authors give a general definition of K\"ahler angle. There are many results about K\"ahler angle one can try to generalize to the general case. In this paper, we focus on the symplectic critical surfaces in Hermite surfaces which is a generalization of \cite{HL1} or \cite{HLS1}. 



\let\thefootnote\relax\footnotetext{\textsuperscript{1} \textit{School of Mathematical Sciences, University of Science and Technology of China, Hefei 230026, P.R.of China}}


\section{Introduction}

Suppose that $(M,J,\langle,\rangle)$ is a Hermite surface, where $\langle,\rangle$ is the Riemannnian metric on $M$, $J$ is a complex structure compatible with $\langle,\rangle$, i.e.$$\langle JU,JV\rangle=\langle U,V\rangle.$$ The K\"ahler form $\omega$ on $M$ is defined by $$\omega(U,V)=\langle JU,V\rangle.$$
For a compact oriented real surface $\Sigma$ which is smoothly immersed in $M$, by the general definition in section $2$, the K\"ahler angle $\alpha$ of $\Sigma$ in $M$ is defined by 
\begin{eqnarray}
   \omega|_{\Sigma}&=&\cos{\alpha}d\mu_{\Sigma},
\end{eqnarray}
 where $d\mu_{\Sigma}$ is the area element of $\Sigma$ of the induced metric $g$ from $<,>$. If $\cos{\alpha}\equiv{1}$, $(\Sigma,J|_{\Sigma},g)$  is a holomorphic curve by the result in section $2$; note that, $\omega|_{\Sigma}$ is the top form, so $d\omega|_{\Sigma}=0$, thus, if $\cos{\alpha}>0$, $\omega|_{\Sigma}$ is the symplectic form, $(\Sigma,J|_{\Sigma},g)$ is symplectic surface; if $\cos{\alpha}\equiv{0}$, $(\Sigma,J|_{\Sigma},g)$ is a Lagrangian surface.

As in \cite{HLS1}, in this paper, under the condition that $M$ is Hermite, we consider a sequence of functionals $$L_{\beta}=\int_{\Sigma}{\frac{1}{\cos^{\beta}{\alpha}}}d\mu.$$

The critical point of the functionals $L_{\beta}$ in the class of symplectic surfaces in a Hermite surface is also called a $\beta-$symplectic critical surface as in \cite{HLS1}. We first calculate the Euler-Lagrange equation of $L_{\beta}$, note that, in our case, $\bar\nabla{J}\ne{0}$, but it is the same as in \cite{HLS1} under a nature condition: 
\begin{equation}\label{Jc1}
  \Xi\ \bar{\nabla}\tilde{J}(X,Y,\xi):=\bar{\nabla}\tilde{J}(X,Y,\xi)+\bar{\nabla}\tilde{J}(Y,\xi,X)+\bar{\nabla}\tilde{J}(\xi,X,Y)=0,
\end{equation}
where $\bar{\nabla}\tilde{J}(X,Y):=\langle\bar{\nabla}J(X),Y\rangle$, $X,Y$ are vector fields of $\Sigma$, $\xi$ is the normal vector field of $\Sigma$ in $M$, $\bar{\nabla}$ is the Levi-Civita connection on $M$. 

\newtheorem{theorem}{Theorem}[section] 

\begin{theorem}
Let $M$ be a Hermite surface. The complex structure $J$ satisfies the condition $(2)$, then the Euler-Lagrange equation of $L_{\beta}(\beta\ne{-1})$ is
\begin{equation}\label{E-L}
  \cos^{3}\alpha\mathbf{H}-\beta(J(J\nabla{\cos\alpha})^{\top})^{\bot} = 0,
\end{equation}
where $H$ is the mean curvature vector of $\Sigma$ in $M$, and $()^{\top}$ means tangential components of $()$, $()^{\bot}$ means the normal components of $()$.
\end{theorem}

Because of the Euler-Lagrange equation is the same as in \cite{HLS1}, so the results which only derived from $(3)$ are believed to hold, such as $(3)$ is an elliptic system modulo tangential diffeomorphisms of $\Sigma$. We then study the harmonic property of $\cos\alpha$ of a symplectic critical surface. Under a nature condition, we have

\begin{theorem}
  Suppose that $M$ is Hermite surface, the complex structure $J$ satisfies $(2)$ and $(\bar\nabla_{X}J(Y)+\bar\nabla_{Y}J(X))^{\bot}=0 $, and $\Sigma$ is a $\beta$-symplectic critical surface in $M$ with the K\"ahker angle $\alpha$. Then $\cos\alpha$ satisfies 
  \begin{eqnarray*}
    \triangle\cos\alpha &=& \frac{2\beta\sin^{2}\alpha}{\cos\alpha(\cos^{2}\alpha+\beta\sin^{2}\alpha)}|\nabla\alpha|^{2}-2\cos\alpha|\nabla\alpha|^{2}  \\
     & & -\frac{\sin\alpha\cos^{2}\alpha}{\cos^{2}\alpha+\beta\sin^{2}\alpha}(K_{1213}-K_{1224}) + \Theta   \\
     & & -\frac{\beta\sin\alpha}{\cos^{2}\alpha+\beta\sin^{2}\alpha}(\partial_{1}\cos\alpha J_{14,2}+3\partial_{2}\cos\alpha J_{13,2}), 
  \end{eqnarray*}
  where
  \begin{eqnarray*}
  \Theta &=& \frac{2\cos\alpha}{\sin^{2}\alpha}(1+\frac{\cos^{2}\alpha}{\cos^{2}\alpha+\beta\sin^{2}\alpha})(\partial_{1}\cos\alpha J_{12,1}+\partial_{2}\cos\alpha J_{12,2})\\
      & & +\frac{\cos^{2}\alpha}{\cos^{2}\alpha+\beta\sin^{2}\alpha}J_{12,kk} - \frac{2\cos^{3}\alpha}{\sin^{2}\alpha(\cos^{2}\alpha+\beta\sin^{2}\alpha)}((J_{12,1})^{2}+(J_{12,2})^{2}).
  \end{eqnarray*}
\end{theorem}


\newpage

\section{The Euler-Lagrange equation of the functional $L_{\beta}(\beta\ne{-1})$ in Hermite surfaces}

Assume that $\phi_{t}:\Sigma\longrightarrow M$ is a one-parameter family of immersions and $\frac{\partial\phi_{t}}{\partial{t}}|_{t=0}=\xi$, $\xi$ is a variational vector field of $\Sigma$. We denote by $\bar{\nabla}$ the Levi-Civita connection and by $K$ the Riemannian curvature tensor on $M$. And $\nabla,R$ denote the Levi-Civita connection and the Riemannian curvature tensor of the induced metric $g$ on $\Sigma$. The proof of the first variational formula is almost the same as that in \cite{HLS1}. We only need to take care of the terms including $\bar\nabla{J}(\ne{0})$.
 
\begin{theorem}
  Let $M$ be a Hermite surface. For any smooth normal vector field $\xi$, if $J$ satisfies the condition $(2)$, the first variational formula of the functional $L_{\beta}$ is
\begin{equation}\label{fvL}
  \delta_{\xi}L_{\beta}= -(\beta+1)\int_{\Sigma}\frac{\xi\cdot\mathbf{H}}{\cos^{\beta}\alpha}d\mu + 
                       \beta(\beta+1)\int_{\Sigma}\frac{\xi\cdot{(J(J\nabla\cos\alpha)^{\top})^{\bot}}}{\cos^{\beta+3}\alpha}d\mu,
\end{equation}
where $\mathbf{H}$ is the mean curvature vector of $\Sigma$ in $M$, and $()^{\top}$ means tangential components of $()$, $()^{\bot}$ means the normal components of $()$. The Euler-Lagrange equation of $L_{\beta}(\beta\ne{-1})$ is
$$\cos^{3}\alpha\mathbf{H}-\beta(J(J\nabla{\cos\alpha})^{\top})^{\bot} = 0.$$
\end{theorem}

\begin{proof}
Let $\{x_{i}\}$ be the local normal coordinates around a fixed point $p$ on $\Sigma$. The induced metric on $\phi_{t}(\Sigma)$ is 
$$g_{ij}(t)=\langle\frac{\partial{\phi_{t}}}{\partial{x_{i}}},\frac{\partial{\phi_{t}}}{\partial{x_{j}}}\rangle.$$
For simplicity, we denote $e_{i}(t)=\frac{\partial{\phi_{t}}}{\partial{x_{i}}}, e_{i}=e_{i}(0)$, and $g_{ij}(t)$ by $g_{ij}$. Then 
$$\frac{\partial}{\partial{t}}|_{t=0}g_{ij}= \frac{\partial}{\partial{t}}|_{t=0}\langle e_{i}(t),e_{j}(t)\rangle
=\langle \bar{\nabla}_{\xi}e_{i},e_{j}\rangle + \langle e_{i},\bar{\nabla}_{\xi}e_{j}\rangle
=\langle \bar{\nabla}_{e_{i}}\xi,e_{j}\rangle + \langle e_{i},\bar{\nabla}_{e_{j}}\xi\rangle$$ 
By the general definition of K\"ahler angle, we have
$$\cos(\alpha_{t})=\frac{\omega(e_{i}(t),e_{j}(t))}{\sqrt{\det(g_{t})}},$$
where $\det(g_{t})$ is the determinant of the metric $(g_{t})$. Then
$$
    \frac{d}{dt}|_{t=0}L_{\beta}(\phi_{t}) =  \int_{\Sigma}\big{(}\frac{1}{\cos^{\beta}{\alpha}}\frac{d}{dt}|_{t=0}(d\mu_{t})-\beta\frac{1}{\cos^{\beta+1}\alpha}\frac{d}{dt}|_{t=0}\cos\alpha d\mu_{t}\big{)}.
$$
A direct calculation gives
$$
\frac{d}{dt}|_{t=0}L_{\beta}(\phi_{t})  =  \int_{\Sigma}\big{(}\frac{(\beta+1)}{2}\frac{1}{\cos^{\beta}{\alpha}}g^{ij}\frac{dg_{ij}}{dt}(0)-\beta\frac{1}{\cos^{\beta+1}\alpha}
\frac{\frac{d}{dt}|_{t=0}\omega(e_{1}(t),e_{2}(t))}{\sqrt{\det(g_{t})}} \big{)}d\mu
$$
By $\xi$ is a normal vector field, we have 
$$
g^{ij}\frac{dg_{ij}}{dt}(0) = g^{ij}(\langle \bar{\nabla}_{e_{i}}\xi,e_{j}\rangle + \langle e_{i},\bar{\nabla}_{e_{j}}\xi\rangle)=2\langle\bar{\nabla}_{e_{i}}\xi,e_{i}\rangle
= -2\langle\bar{\nabla}_{e_{i}}e_{i},\xi\rangle=-2\xi\cdot\mathbf{H},
$$
\begin{eqnarray}\label{w}
  \frac{d}{dt}|_{t=0}\omega(e_{1}(t),e_{2}(t)) & = & \bar{\nabla}_{\xi}\omega(e_{1},e_{2}) + \omega(\bar{\nabla}_{\xi}e_{1},e_{2}) + \omega(e_{1},\bar{\nabla}_{\xi}e_{2}) \nonumber \\
   & = & \bar{\nabla}_{\xi}\omega(e_{1},e_{2}) + \omega(\bar{\nabla}_{e_{1}}\xi,e_{2}) + \omega(e_{1},\bar{\nabla}_{e_{2}}\xi).
\end{eqnarray}
So, 
$$
\begin{array}{ccc}
\frac{d}{dt}|_{t=0}L_{\beta}(\phi_{t}) & = &-(\beta+1)\int_{\Sigma}\frac{\xi\cdot\mathbf{H}}{\cos^{\beta}{\alpha}}d\mu \\
& & -\beta\int_{\Sigma}\frac{\bar{\nabla}_{\xi}\omega(e_{1},e_{2}) + \omega(\bar{\nabla}_{e_{1}}\xi,e_{2}) + \omega(e_{1},\bar{\nabla}_{e_{2}}\xi)}{\cos^{\beta+1}\alpha}dx^{1}\wedge dx^{2}
\end{array}
$$
By the condition $(\ref{Jc1})$, $\bar{\nabla}\omega(U,V)=\langle \bar{\nabla}J(U),V\rangle = \bar{\nabla}\tilde{J}(U,V)$, and
$$
\omega(\bar{\nabla}_{e_{1}}\xi,e_{2})=e_{1}(\omega(\xi,e_{2}))-\bar{\nabla}_{e_{1}}\omega(\xi,e_{2})-\omega(\xi,\bar{\nabla}_{e_{1}}e_{2}),
$$
$$
\omega(e_{1},\bar{\nabla}_{e_{2}}\xi)=e_{2}(\omega(e_{1},\xi))-\bar{\nabla}_{e_{2}}\omega(e_{1},\xi)-\omega(\bar{\nabla}_{e_{2}}e_{1},\xi),
$$
we have 
\begin{eqnarray*}
  (\ref{w}) & = & \bar{\nabla}_{\xi}\omega(e_{1},e_{2}) - \bar{\nabla}_{e_{1}}\omega(\xi,e_{2}) - \bar{\nabla}_{e_{2}}\omega(e_{1},\xi)\\
      & & + e_{1}(\omega(\xi,e_{2})) + e_{2}(\omega(e_{1},\xi))- \omega(\xi,\bar{\nabla}_{e_{1}}e_{2}-\bar{\nabla}_{e_{2}}e_{1}) \\
      & = & \bar{\nabla}_{\xi}\omega(e_{1},e_{2}) + \bar{\nabla}_{e_{1}}\omega(e_{2},\xi) + \bar{\nabla}_{e_{2}}\omega(\xi,e_{1})\\
      & & + e_{1}(\omega(\xi,e_{2})) + e_{2}(\omega(e_{1},\xi))\\
      & = & \bar{\nabla}\omega(e_{1},e_{2},\xi) + \bar{\nabla}\omega(e_{2},\xi,e_{1}) + \bar{\nabla}\omega(\xi,e_{1},e_{2})\\
      & & + e_{1}(\omega(\xi,e_{2})) + e_{2}(\omega(e_{1},\xi))\\
      & = & \Im\ {\bar{\nabla}\tilde{J}}(e_{1},e_{2},\xi) + e_{1}(\omega(\xi,e_{2})) + e_{2}(\omega(e_{1},\xi))\\
      & = & e_{1}(\omega(\xi,e_{2})) + e_{2}(\omega(e_{1},\xi)).
\end{eqnarray*}
Since $\Sigma$ is closed, applying the Stokes formula, we obtain 
\begin{eqnarray*}
& & -\beta\int_{\Sigma}\frac{\bar{\nabla}_{\xi}\omega(e_{1},e_{2}) + \omega(\bar{\nabla}_{e_{1}}\xi,e_{2}) + \omega(e_{1},\bar{\nabla}_{e_{2}}\xi)}{\cos^{\beta+1}\alpha}dx^{1}\wedge dx^{2}\\
&=& -\beta\int_{\Sigma}\frac{e_{1}(\omega(\xi,e_{2})) + e_{2}(\omega(e_{1},\xi))}{\cos^{\beta+1}\alpha}dx^{1}\wedge dx^{2}\\
&=& -\beta(\beta+1)\int_{\Sigma}\frac{\omega(\xi,e_{2})\nabla_{e_{1}}\cos\alpha + \omega(e_{1},\xi)\nabla_{e_{2}}\cos\alpha}{\cos^{\beta+2}\alpha}dx^{1}\wedge dx^{2}.
\end{eqnarray*}
Since $\omega(\xi,e_{2})=-\langle \xi,Je_{2}\rangle$, $\omega(e_{1},\xi)=\langle \xi,Je_{1}\rangle$, $\cos\alpha=\langle Je_{1},e_{2}\rangle=-\langle e_{1},Je_{2}\rangle$, we have
\begin{equation}\label{int}
\omega(\xi,e_{2})\nabla_{e_{1}}\cos\alpha + \omega(e_{1},\xi)\nabla_{e_{2}}\cos\alpha = -\langle \xi, J(e_{2}\nabla_{e_{1}}\cos\alpha-e_{1}\nabla_{e_{2}}\cos\alpha)\rangle
\end{equation}
and 
\begin{eqnarray*}
  (J\nabla\cos\alpha)^{\top} &=& (Je_{1}\nabla_{e_{1}}\cos\alpha+Je_{2}\nabla_{e_{2}}\cos\alpha)^{\top} \\
    &=& \langle Je_{1},e_{2}\rangle e_{2}\nabla_{e_{1}}\cos\alpha + \langle Je_{2},e_{1}\rangle e_{1}\nabla_{e_{2}}\cos\alpha\\
    &=& (e_{2}\nabla_{e_{1}}\cos\alpha-e_{1}\nabla_{e_{2}}\cos\alpha)\cos\alpha.
\end{eqnarray*}
So 
$$(\ref{int})= -\langle \xi,J(J\nabla\cos\alpha)^{\top} \rangle/ \cos\alpha = -\frac{\xi\cdot(J(J\nabla\cos\alpha)^{\top})^{\bot}}{\cos\alpha}.$$
Therefore, we have 
\begin{eqnarray*}
  \frac{d}{dt}|_{t=0}L_{\beta}(\phi_{t}) &=&  -(\beta+1)\int_{\Sigma}\frac{\xi\cdot\mathbf{H}}{\cos^{\beta}\alpha}d\mu \\
  & & +  \beta(\beta+1)\int_{\Sigma}\frac{\xi\cdot{(J(J\nabla\cos\alpha)^{\top})^{\bot}}}{\cos^{\beta+3}\alpha}d\mu.
\end{eqnarray*}

\end{proof}

We express $(J(J\nabla\cos\alpha)^{\top})^{\bot}$ at a fixed point p in a local frame. Note that $\omega(e_{i},e_{j})=\langle Je_{i},e_{j}\rangle=:J_{ij}$, $J_{ji}=\omega(e_{j},e_{i})=-\omega(e_{i},e_{j})=-J_{ij}$, so $J^{T}=-J$, $JJ^{T}=-J^{2}=id$. Then a simple calculation show that $J$ has the form
$$\left(
\begin{array}{cccc}
   0 & x & y & z \\
  -x & 0 & z & -y \\
  -y & -z & 0 & x \\
  -z & y & -x & 0 
\end{array}\right),
$$
where $x^{2}+y^{2}+z^{2}=1$. By the definition of the K\"ahler angle, we know that at p, $$x=\langle Je_{1},e_{2}\rangle=\omega(e_{1},e_{2})=\cos\alpha.$$
Then 
\begin{eqnarray*}
  (J(J\nabla\cos\alpha)^{\top})^{\bot} &=& (J(\cos\alpha\partial_{1}\cos\alpha e_{2} - \cos\alpha\partial_{2}\cos\alpha e_{1}))^{\bot} \\
   &=&  \cos\alpha\partial_{1}\cos\alpha (Je_{2})^{\bot} - \cos\alpha\partial_{2}\cos\alpha (Je_{1})^{\bot}\\
   &=&  \cos\alpha\partial_{1}\cos\alpha (ze_{3}-ye_{4}) - \cos\alpha\partial_{2}\cos\alpha (ye_{3}+ze_{4}).
\end{eqnarray*}
Thus, by $(\ref{E-L})$ we get that 
\begin{eqnarray}
  H^{3} &=& -\beta\frac{1}{\cos^{2}\alpha}(y\partial_{2}\cos\alpha - z\partial_{1}\cos\alpha); \\
  H^{4} &=& -\beta\frac{1}{\cos^{2}\alpha}(y\partial_{1}\cos\alpha + z\partial_{2}\cos\alpha).
\end{eqnarray}

Because of the Euler-Lagrange equation is the same as in \cite{HLS1}, one believe that $(\ref{E-L})$ is an elliptic system modulo tangential diffeomorphisms of $\Sigma$. In fact, if one check the proof in the section $3$ of \cite{HLS1}, the only difference is the additional first-order terms caused by $\bar\nabla J\ne{0}$ which have no effect on the calculation of the principal symbol.

\newtheorem{proposition}[theorem]{Proposition}

\begin{proposition}
For a $\beta$-symplectic critical surface $\Sigma$ in a Hermite surface $M$ with $\beta\ge0$, the Euler-Lagrange equation $(\ref{E-L})$ is an elliptic system modulo tangential diffeomorphisms of $\Sigma$.
\end{proposition}


\section{The harmonic property of $\cos\alpha$ of a symplectic critical surface}

Let $T\Sigma$, $N\Sigma$ be the tangent bundle and the normal bundle of $\Sigma$ in $M$ respectively. The second fundamental form $B:T\Sigma\times T\Sigma\to N\Sigma$ is defined by $B(X,Y)=(\bar{\nabla}_{X}Y)^{\bot}$ for any tangent vector fields $X,Y$. The operator $A:T\Sigma\times N\Sigma\to T\Sigma$ is defined by $A_{\xi}(X)=(\bar\nabla_{X}\xi)^{\top}$, $\xi\in N\Sigma$. Here $()^{\bot}$ denotes the projection from $TM$ onto $N\Sigma$, $()^{\top}$ denotes the projection onto $T\Sigma$. Evidently, $$\langle B(X,Y),\xi\rangle = -\langle A_{\xi}(X),Y\rangle.$$
The Gauss formula is 
$$\bar{\nabla}_{X}Y = \nabla_{X}Y + B(X,Y),$$
where $\nabla$ is the Levi-Civita connection on $\Sigma$. 
The Weingarten formula is
$$\bar{\nabla}_{X}\xi = -A_{\xi}X + \nabla^{\bot}_{X}\xi,$$
where $\nabla^{\bot}$ is the normal connection of $N\Sigma$. 

\begin{proposition}\label{propxint}
Let $(M,J,\langle,\rangle)$ is a Hermite surface, where $\langle,\rangle$ is the Riemannnian metric on $M$, $J$ is a complex structure compatible with $\langle,\rangle$. If $\Sigma$ is a closed symplectic surface which is smoothly immersed in $M$ with the K\"ahler angle $\alpha$, then 
\begin{eqnarray}\label{trieq}
  \triangle\cos\alpha &=& -\cos\alpha(|h^{3}_{1k}-h^{4}_{2k}|^{2}+|h^{4}_{1k}+h^{3}_{2k}|^{2})+ \sin\alpha(H^{4}_{,1}+H^{3}_{,2})\nonumber\\
   & & -\sin\alpha(K_{1213}-K_{1224}) + J_{12,kk}+2J_{\alpha 2,k}h^{\alpha}_{1k}+2J_{1\alpha,k}h^{\alpha}_{2k},
\end{eqnarray}
where $K$ is the curvature tensor of $M$ and $H^{\alpha}_{,i}=\langle \bar\nabla_{e_{i}}H,e_{\alpha}\rangle$.
\end{proposition}

\begin{proof}
For a fixed point $p\in\Sigma$, we first choose $\{x_{i}\}$ be the local normal coordinates around $p\in\Sigma$, the frame $\{\partial_{i}\}$ are denote by $\{e_{i}\}, i=1,2$, so $\nabla_{e_{i}}e_{j}(p)=0$. Then we choose normal frame $\{e_{\alpha}\},\alpha=3,4$ such that at $p$, $y=\sin\alpha,z=0$ and $\nabla^{\bot}_{e_{i}}e_{\alpha}(p)=0$. 

By the definition of K\"ahler angle, 
\begin{eqnarray}\label{grad}
\nabla_{e_{k}}\cos\alpha=\nabla_{e_{k}}(\omega(e_{1},e_{2}))/\sqrt{\det g}-\frac{1}{2}\cos\alpha g^{ij}\nabla_{e_{k}}(g_{ij}).
\end{eqnarray}
Then, at $p$, we have 
\begin{eqnarray}\label{tri}
  \triangle\cos\alpha &=& \nabla_{e_{k}}(\nabla_{e_{k}}\cos\alpha) \nonumber\\
   &=& \nabla_{e_{k}}(\nabla_{e_{k}}(\omega(e_{1},e_{2}))/\sqrt{\det g}-\frac{1}{2}\cos\alpha g^{ij}\nabla_{e_{k}}(g_{ij})) \nonumber\\
   &=& \triangle(\omega(e_{1},e_{2}))-\frac{1}{2}\cos\alpha\triangle(g_{ij})g^{ij}.
\end{eqnarray}
Note that $\bar\nabla\omega\ne{0}$, we have 
$$\bar{\nabla}_{e_{k}}(\omega(e_{1},e_{2}))=\bar{\nabla}_{e_{k}}\omega(e_{1},e_{2})+\omega(\bar{\nabla}_{e_{k}}e_{1},e_{2})+\omega(e_{1},\bar{\nabla}_{e_{k}}e_{2}),$$
$$\bar{\nabla}_{e_{k}}\omega(e_{1},e_{2})=\langle \bar{\nabla}_{e_{k}}J(e_{1}),e_{2}\rangle.$$
Then, we have 
\begin{eqnarray*}
  \triangle(\omega(e_{1},e_{2}))&=& \bar{\nabla}_{e_{k}}(\bar{\nabla}_{e_{k}}(\omega(e_{1},e_{2}))) \\
   &=& \langle\bar{\nabla}_{e_{k}}\bar{\nabla}_{e_{k}}J(e_{1}),e_{2}\rangle + \langle\bar{\nabla}_{e_{k}}J(\bar{\nabla}_{e_{k}}e_{1}),e_{2}\rangle+
       \langle\bar{\nabla}_{e_{k}}J(e_{1}),\bar{\nabla}_{e_{k}}e_{2}\rangle  \\
   & & + \bar{\nabla}_{e_{k}}\omega(\bar{\nabla}_{e_{k}}e_{1},e_{2})+ \omega(\bar{\nabla}_{e_{k}}\bar{\nabla}_{e_{k}}e_{1},e_{2}) + \omega(\bar{\nabla}_{e_{k}}e_{1},\bar{\nabla}_{e_{k}}e_{2})  \\
   & & + \bar{\nabla}_{e_{k}}\omega(e_{1},\bar{\nabla}_{e_{k}}e_{2}) + \omega(\bar{\nabla}_{e_{k}}e_{1},\bar{\nabla}_{e_{k}}e_{2}) + 
   \omega(e_{1},\bar{\nabla}_{e_{k}}\bar{\nabla}_{e_{k}}e_{2}) \\
   &=& J_{12,kk} + 2J_{\alpha2,k}h^{\alpha}_{1k} + 2J_{1\alpha,k}h^{\alpha}_{2k} \\
   & & + \omega(\bar{\nabla}_{e_{k}}\bar{\nabla}_{e_{k}}e_{1},e_{2}) + \omega(e_{1},\bar{\nabla}_{e_{k}}\bar{\nabla}_{e_{k}}e_{2}) + 2\omega(\bar{\nabla}_{e_{k}}e_{1},\bar{\nabla}_{e_{k}}e_{2}) \\
   &=:& J_{12,kk} + 2J_{\alpha2,k}h^{\alpha}_{1k} + 2J_{1\alpha,k}h^{\alpha}_{2k} + (*1),
\end{eqnarray*}
and 
\begin{eqnarray*}
  (*1) &=& \omega(\bar{\nabla}_{e_{k}}({\nabla}_{e_{k}}e_{1} + B(e_{k},e_{1})),e_{2}) + \omega(e_{1},\bar{\nabla}_{e_{k}}({\nabla}_{e_{k}}e_{2}+B(e_{k},e_{2})))  \\
   & & + 2\omega(B(e_{k},e_{1}),B({e_{k}},e_{2})) \\
   &=& \omega(\nabla_{e_{k}}{\nabla}_{e_{k}}e_{1},e_{2}) + \omega(B(e_{k},\nabla_{e_{k}}e_{1}),e_{2}) + \omega(\bar\nabla_{e_{k}}(B(e_{k},e_{1})),e_{2}) \\
   & & \omega(e_{1},\nabla_{e_{k}}{\nabla}_{e_{k}}e_{2}) + \omega(e_{1},B(e_{k},\nabla_{e_{k}}e_{2})) + \omega(e_{1},\bar\nabla_{e_{k}}(B(e_{k},e_{2}))) \\
   & & + 2\omega(B(e_{k},e_{1}),B({e_{k}},e_{2})) \\
   &=&  \cos\alpha\langle \nabla_{e_{k}}{\nabla}_{e_{k}}e_{1},e_{1}\rangle + \cos\alpha\langle \nabla_{e_{k}}{\nabla}_{e_{k}}e_{2},e_{2}\rangle \\
   & & + \omega(\bar\nabla_{e_{k}}(B(e_{k},e_{1})),e_{2}) + \omega(e_{1},\bar\nabla_{e_{k}}(B(e_{k},e_{2}))) + 2\omega(B(e_{k},e_{1}),B({e_{k}},e_{2}))\\
   &=:&  \cos\alpha\langle \nabla_{e_{k}}{\nabla}_{e_{k}}e_{1},e_{1}\rangle + \cos\alpha\langle \nabla_{e_{k}}{\nabla}_{e_{k}}e_{2},e_{2}\rangle +(*2),
\end{eqnarray*}
and 
\begin{eqnarray*}
  (*2) &=& \omega(\bar\nabla_{e_{k}}(h^{\alpha}_{1k}e_{\alpha}),e_{2}) + \omega(e_{1},\bar\nabla_{e_{k}}(h^{\alpha}_{2k}e_{\alpha})) + 2\omega(h^{\alpha}_{1k}e_{\alpha},h^{\beta}_{2k}e_{\beta}) \\
   &=& \omega(h^{\alpha}_{k1,k}e_{\alpha}-h^{\alpha}_{1k}h^{\alpha}_{lk}e_{l},e_{2}) - \omega(h^{\alpha}_{k2,k}e_{\alpha}-h^{\alpha}_{2k}h^{\alpha}_{lk}e_{l},e_{1}) \\
   & & + 2\omega(h^{\alpha}_{1k}e_{\alpha},h^{\beta}_{2k}e_{\beta})\\
   &=& \cos\alpha(-(h^{\alpha}_{1k})^{2}-(h^{\alpha}_{2k})^{2} + 2h^{3}_{1k}h^{4}_{2k} - 2h^{4}_{1k}h^{3}_{2k}) \\
   & & + \omega(h^{\alpha}_{kk,1}e_{\alpha}-K_{\alpha k1k}e_{\alpha},e_{2})-\omega(h^{\alpha}_{kk,2}e_{\alpha}-K_{\alpha k2k}e_{\alpha},e_{1})\\
   &=&  \cos\alpha(-(h^{\alpha}_{1k})^{2}-(h^{\alpha}_{2k})^{2} + 2h^{3}_{1k}h^{4}_{2k} - 2h^{4}_{1k}h^{3}_{2k}) \\
   & & + \sin\alpha(H^{4}_{,1}+H^{3}_{,2}) - \sin\alpha(K_{4k1k}+K_{3k2k}) \\
   &=&  \cos\alpha(-(h^{\alpha}_{1k})^{2}-(h^{\alpha}_{2k})^{2} + 2h^{3}_{1k}h^{4}_{2k} - 2h^{4}_{1k}h^{3}_{2k}) \\
   & & + \sin\alpha(H^{4}_{,1}+H^{3}_{,2}) - \sin\alpha(K_{1213}-K_{1224}).
\end{eqnarray*}
It is easy to check that 
\begin{eqnarray*}
  \frac{1}{2}\cos\alpha\triangle(g_{ij})g^{ij} &=& \frac{1}{2}\cos\alpha\triangle \langle e_{i},e_{j}\rangle g^{ij} \\
   &=& \cos\alpha\langle \nabla_{e_{k}}\nabla_{e_{k}}e_{i},e_{j}\rangle g^{ij} \\
   &=& \cos\alpha\langle \nabla_{e_{k}}\nabla_{e_{k}}e_{1},e_{1}\rangle + \cos\alpha\langle \nabla_{e_{k}}\nabla_{e_{k}}e_{2},e_{2}\rangle.
\end{eqnarray*}
By $(\ref{tri})$, we obtain, 
\begin{eqnarray*}
  \triangle\cos\alpha &=& \cos\alpha(-(h^{\alpha}_{1k})^{2}-(h^{\alpha}_{2k})^{2} + 2h^{3}_{1k}h^{4}_{2k} - 2h^{4}_{1k}h^{3}_{2k}) \\
   & &  + \sin\alpha(H^{4}_{,1}+H^{3}_{,2}) - \sin\alpha(K_{1213}-K_{1224})\\
   & &  + J_{12,kk} + 2J_{\alpha2,k}h^{\alpha}_{1k} + 2J_{1\alpha,k}h^{\alpha}_{2k}.
\end{eqnarray*}

\end{proof}

By $(\ref{grad})$, at $p$, we have 
\begin{eqnarray}\label{partial1}
  \partial_{1}\cos\alpha &=& \bar\nabla_{e_{1}}\cos\alpha =\nabla_{e_{1}}\cos\alpha = \bar\nabla_{e_{1}}(\omega(e_{1},e_{2})) \nonumber \\
   &=& \bar\nabla_{e_{1}}\omega(e_{1},e_{2}) + \omega(\bar\nabla_{e_{1}}e_{1},e_{2}) +\omega(e_{1},\bar\nabla_{e_{1}}e_{2}) \nonumber \\
   &=& J_{12,1} - \langle \bar\nabla_{e_{1}}e_{1},Je_{2}\rangle + \langle Je_{1}, \bar\nabla_{e_{1}}e_{2}\rangle \nonumber \\
   &=& J_{12,1} +\sin\alpha (h_{11}^{4}+h_{12}^{3}).
\end{eqnarray}
Similarly,
\begin{eqnarray}\label{partial2}
\partial_{2}\cos\alpha =  J_{12,2} +\sin\alpha (h_{12}^{4}+h_{22}^{3}).
\end{eqnarray}
When we assume $\Sigma$ is a $\beta$-symplectic critical surface, the Euler-Lagrange equation or the last two identities show that $h^{\alpha}_{ij}$ is connected with $\partial_{i}\cos\alpha$. We first 
consider $J_{\alpha 2,k}h^{\alpha}_{1k}+J_{1\alpha,k}h^{\alpha}_{2k}$ in $(\ref{trieq})$:
\begin{eqnarray*}
  J_{\alpha 2,k}h^{\alpha}_{1k}+J_{1\alpha,k}h^{\alpha}_{2k} &=& 
       J_{32,1}h^{3}_{11}+J_{32,2}h^{3}_{12}+J_{42,1}h^{4}_{11}+J_{42,2}h^{4}_{12} \\
   & & +J_{13,1}h^{3}_{21}+J_{13,2}h^{3}_{22}+J_{14,1}h^{4}_{21}+J_{14,2}h^{4}_{22} \\
   &=& J_{32,1}h^{3}_{11}+J_{42,1}h^{4}_{11}+J_{13,2}h^{3}_{22}+J_{14,2}h^{4}_{22} \\
   & & +(J_{13,1}+J_{32,2})h^{3}_{12}+(J_{14,1}+J_{42,2})h^{4}_{12} \\
   &=& J_{32,1}H^{3}+J_{14,2}H^{4} \\
   & & +(J_{13,2}-J_{32,1})h^{3}_{22}+(J_{42,1}-J_{14,2})h^{4}_{11}\\
   & & +(J_{13,1}+J_{32,2})h^{3}_{12}+(J_{14,1}+J_{42,2})h^{4}_{12} \\
   &=& J_{32,1}H^{3}+J_{14,2}H^{4} \\
   & & -(J_{13,2}-J_{32,1})(\frac{\partial_{2}\cos\alpha}{-\sin\alpha}+h_{12}^{4}+\frac{J_{12,2}}{\sin\alpha})\\
   & & -(J_{42,1}-J_{14,2})(\frac{\partial_{1}\cos\alpha}{-\sin\alpha}+h_{12}^{3}+\frac{J_{12,1}}{\sin\alpha})\\
   & & +(J_{13,1}+J_{32,2})h^{3}_{12}+(J_{14,1}+J_{42,2})h^{4}_{12} \\
   &=& J_{32,1}H^{3}+J_{14,2}H^{4} \\
   & & -(J_{13,2}-J_{32,1})(\frac{\partial_{2}\cos\alpha}{-\sin\alpha}+\frac{J_{12,2}}{\sin\alpha})\\
   & & -(J_{42,1}-J_{14,2})(\frac{\partial_{1}\cos\alpha}{-\sin\alpha}+\frac{J_{12,1}}{\sin\alpha})\\
   & & +(J_{13,1}+J_{32,2}-J_{42,1}+J_{14,2})h^{3}_{12}\\
   & & +(J_{14,1}+J_{42,2}-J_{13,2}+J_{32,1})h^{4}_{12}.
\end{eqnarray*}
So we need 
\begin{eqnarray*}
  J_{13,1}+J_{32,2}-J_{42,1}+J_{14,2} &=& J^{3}_{1,1}-J^{3}_{2,2}+J^{4}_{2,1}+J^{4}_{1,2} = 0 \\
  J_{14,1}+J_{42,2}-J_{13,2}+J_{32,1} &=& J^{4}_{1,1}-J^{4}_{2,2}-(J^{3}_{1,2}+J^{3}_{2,1}) = 0.
\end{eqnarray*}
A natural condition is $J^{\alpha}_{i,j}+J^{\alpha}_{j,i}=0$, which gives that $J^{\alpha}_{i,i}=0$. That is 
\begin{eqnarray}\label{Jc2}
 (\bar\nabla_{X}J(Y)+\bar\nabla_{Y}J(X))^{\bot}=0
\end{eqnarray} 
Moreover, $J_{13,2}=J^{3}_{1,2}=-J^{3}_{2,1}=-J_{23,1}=J_{32,1}$, similarly, $J_{42,1}=J_{14,2}$. Thus
\begin{eqnarray*}
  J_{\alpha 2,k}h^{\alpha}_{1k}+J_{1\alpha,k}h^{\alpha}_{2k} &=& J_{13,2}H^{3}+J_{14,2}H^{4} \\
    &=& -\beta\frac{\sin\alpha}{\cos^{2}\alpha}(J_{13,2}\partial_{2}\cos\alpha + J_{14,2}\partial_{1}\cos\alpha).
\end{eqnarray*}
Note that the conditions, i.e. $(\ref{Jc1})$ and $(\ref{Jc2})$, can not give $\bar\nabla J=0$.

\begin{theorem}
  Suppose that $M$ is Hermite surface, the complex structure $J$ satisfies $(\ref{Jc1})$ and $(\ref{Jc2})$, and $\Sigma$ is a $\beta$-symplectic critical surface in $M$ with the K\"ahker angle $\alpha$. Then $\cos\alpha$ satisfies 
  \begin{eqnarray}\label{strieq}
    \triangle\cos\alpha &=& \frac{2\beta\sin^{2}\alpha}{\cos\alpha(\cos^{2}\alpha+\beta\sin^{2}\alpha)}|\nabla\alpha|^{2}-2\cos\alpha|\nabla\alpha|^{2} \nonumber \\
     & & -\frac{\sin\alpha\cos^{2}\alpha}{\cos^{2}\alpha+\beta\sin^{2}\alpha}(K_{1213}-K_{1224}) + \Theta  \nonumber \\
     & & -\frac{\beta\sin\alpha}{\cos^{2}\alpha+\beta\sin^{2}\alpha}(\partial_{1}\cos\alpha J_{14,2}+3\partial_{2}\cos\alpha J_{13,2}), 
  \end{eqnarray}
  where
  \begin{eqnarray*}
  \Theta &=& \frac{2\cos\alpha}{\sin^{2}\alpha}(1+\frac{\cos^{2}\alpha}{\cos^{2}\alpha+\beta\sin^{2}\alpha})(\partial_{1}\cos\alpha J_{12,1}+\partial_{2}\cos\alpha J_{12,2})\\
      & & +\frac{\cos^{2}\alpha}{\cos^{2}\alpha+\beta\sin^{2}\alpha}J_{12,kk} - \frac{2\cos^{3}\alpha}{\sin^{2}\alpha(\cos^{2}\alpha+\beta\sin^{2}\alpha)}((J_{12,1})^{2}+(J_{12,2})^{2}).
  \end{eqnarray*}
\end{theorem}

\begin{proof}
We choose the local frame as in Proposition $\ref{propxint}$. For a $\beta$-symplectic critical surface $\Sigma$, if we set
$\tilde{V}=\nabla_{e_{2}}\cos\alpha e_{3}+\nabla_{e_{1}}\cos\alpha e_{4}:=V+\delta V$, where $V=\sin\alpha((h^{4}_{12}+h^{3}_{22})e_{3}+(h^{4}_{11}+h^{3}_{12})e_{4})$,
$\delta V=J_{12,2}e_{3}+J_{12,1}e_{4}$.
Then we have 
$$H=-\beta\frac{\sin\alpha}{\cos^{2}\alpha}\tilde V.$$
By a direct computation, we have at $p$, 
\begin{eqnarray*}
  |h^{3}_{1k}-h^{4}_{2k}|^{2}+|h^{4}_{1k}+h^{3}_{2k}|^{2} &=& |H-\frac{V}{\sin\alpha}|^{2} + |\frac{V}{\sin\alpha}|^{2} \\
   &=& |(\beta\frac{\sin^{2}\alpha}{\cos^{2}\alpha}+1)\frac{\tilde V}{\sin\alpha}-\frac{\delta V}{\sin\alpha}|^{2} + |\frac{\tilde V}{\sin\alpha} - \frac{\delta V}{\sin\alpha}|^{2} \\
   &=& \frac{\beta^{2}\sin^{4}\alpha+2\cos^{4}\alpha+2\beta\sin^{\alpha}\cos^{2}\alpha}{\cos^{4}\alpha}|\nabla\alpha|^{2} \\
   & & +\frac{2}{\sin^{2}\alpha}((J_{12,1})^{2}+(J_{12,2})^{2}) \\
   & & -\frac{2(2\cos^{2}\alpha+\beta\sin^{2}\alpha)}{\cos^{2}\alpha\sin^{2}\alpha}(J_{12,1}\partial_{1}\cos\alpha + J_{12,2}\partial_{2}\cos\alpha),
\end{eqnarray*}
and 
\begin{eqnarray*}
  \sin\alpha(H_{,1}^{4}+H_{,2}^{3}) &=& \sin\alpha(\langle \bar\nabla_{e_{1}}H, e_{4}\rangle+\langle\bar\nabla_{e_{2}}H, e_{3}\rangle)  \\
   &=& \sin\alpha(\partial_{1}H^{4} + H^{3}\langle \bar\nabla_{e_{1}}e_{3},e_{4}\rangle + \partial_{2}H^{3} + H^{4}\langle \bar\nabla_{e_{2}}e_{4},e_{3}\rangle )  \\
   &=& -\beta\sin\alpha(\partial_{1}(\frac{1}{\cos^{2}\alpha}(y\partial_{1}\cos\alpha+z\partial_{2}\cos\alpha))  \\ 
   & & +\partial_{2}(\frac{1}{\cos^{2}\alpha}(y\partial_{2}\cos\alpha-z\partial_{1}\cos\alpha)))   \\
   &=& -\beta\sin\alpha(\partial_{1}(\frac{\partial_{1}\cos\alpha}{\cos^{2}\alpha})\sin\alpha + \partial_{1}y\frac{\partial_{1}\cos\alpha}{\cos^{2}\alpha} +\partial_{1}z\frac{\partial_{2}\cos\alpha}{\cos^{2}\alpha}   \\
   & & + \partial_{2}(\frac{\partial_{2}\cos\alpha}{\cos^{2}\alpha})\sin\alpha + \partial_{2}y\frac{\partial_{2}\cos\alpha}{\cos^{2}\alpha} -\partial_{2}z\frac{\partial_{1}\cos\alpha}{\cos^{2}\alpha})   \\   
   &=& -\beta\sin^{2}\alpha(\partial_{1}(\frac{\partial_{1}\cos\alpha}{\cos^{2}\alpha}) + \partial_{2}(\frac{\partial_{2}\cos\alpha}{\cos^{2}\alpha}))  \\
   & & -\beta\frac{\sin\alpha}{\cos^{2}\alpha}(\partial_{1}y\partial_{1}\cos\alpha+\partial_{2}y\partial_{2}\cos\alpha \\
   & & +\partial_{1}z\partial_{2}\cos\alpha - \partial_{2}z\partial_{1}\cos\alpha) \\
   &=& -\beta\frac{\sin^{2}\alpha}{\cos^{2}\alpha}\triangle\cos\alpha + \frac{2\beta\sin^{2}\alpha}{\cos^{3}\alpha}|\nabla\cos\alpha|^{2} \\
   & & -\beta\frac{\sin\alpha}{\cos^{2}\alpha}(\partial_{1}y\partial_{1}\cos\alpha+\partial_{2}y\partial_{2}\cos\alpha \\
   & & +\partial_{1}z\partial_{2}\cos\alpha - \partial_{2}z\partial_{1}\cos\alpha). \\
\end{eqnarray*}

Now we begin to compute $\partial_{1}y\partial_{1}\cos\alpha+\partial_{2}y\partial_{2}\cos\alpha $ and $\partial_{1}z\partial_{2}\cos\alpha - \partial_{2}z\partial_{1}\cos\alpha$. Note that $y=\langle Je_{1},e_{3}\rangle$, $z=\langle Je_{1},e_{4}\rangle$,$J^{\alpha}_{i,i}=0$. Then, by $(\ref{partial1}),(\ref{partial2})$,
\begin{eqnarray*}
  \partial_{1}y\partial_{1}\cos\alpha+\partial_{2}y\partial_{2}\cos\alpha &=& \partial_{1}\cos\alpha(\bar\nabla_{e_{1}}\langle Je_{1},e_{3}\rangle) +  \partial_{2}\cos\alpha(\bar\nabla_{e_{2}}\langle Je_{1},e_{3}\rangle)   \\
   &=& \partial_{1}\cos\alpha(J_{13,1} + \langle J\bar\nabla_{e_{1}}e_{1},e_{3}\rangle + \langle Je_{1},\bar\nabla_{e_{1}}e_{3}\rangle)    \\
   & & + \partial_{2}\cos\alpha(J_{13,2} + \langle J\bar\nabla_{e_{2}}e_{1},e_{3}\rangle + \langle Je_{1},\bar\nabla_{e_{2}}e_{3}\rangle)  \\
   &=& \partial_{1}\cos\alpha(-\cos\alpha h^{4}_{11}-\cos\alpha h_{11}^{3})    \\
   & & + \partial_{2}\cos\alpha(J_{13,2} -\cos\alpha h^{4}_{12}-\cos\alpha h_{22}^{3})  \\
   &=& \cos\alpha\partial_{1}\cos\alpha(-h^{4}_{11}-h_{11}^{3})    \\
   & & + \cos\alpha\partial_{2}\cos\alpha(- h^{4}_{12} - h_{22}^{3}) + \partial_{2}\cos\alpha J_{13,2}  \\
   &=& \cos\alpha\partial_{1}\cos\alpha(-\frac{\partial_{1}\cos\alpha}{\sin\alpha}+\frac{J_{12,1}}{\sin\alpha})    \\
   & & + \cos\alpha\partial_{2}\cos\alpha(-\frac{\partial_{2}\cos\alpha}{\sin\alpha}+\frac{J_{12,2}}{\sin\alpha}) \\
   & & + \partial_{2}\cos\alpha J_{13,2}  \\
   &=& -\frac{\cos\alpha}{\sin\alpha}|\nabla\cos\alpha|^{2}+\partial_{1}\cos\alpha\frac{J_{12,1}}{\tan\alpha}    \\
   & & + \partial_{2}\cos\alpha(\frac{J_{12,2}}{\tan\alpha} + J_{13,2}).  
\end{eqnarray*}
  
Similarly, 
\begin{eqnarray*}
  \partial_{1}z\partial_{2}\cos\alpha - \partial_{2}z\partial_{1}\cos\alpha &=&  \partial_{2}\cos\alpha(\bar\nabla_{e_{1}}\langle Je_{1},e_{4}\rangle) -  \partial_{1}\cos\alpha(\bar\nabla_{e_{2}}\langle Je_{1},e_{4}\rangle)   \\
   &=& \partial_{2}\cos\alpha(J_{14,1}+\langle J\bar\nabla_{e_{1}}e_{1},e_{4}\rangle +\langle Je_{1},\bar\nabla_{e_{1}}e_{4}\rangle)   \\
   & & - \partial_{1}\cos\alpha(J_{14,2} +\langle J\bar\nabla_{e_{2}}e_{1},e_{4}\rangle + \langle Je_{1},\bar\nabla_{e_{2}}e_{4}\rangle)  \\
   &=& \partial_{2}\cos\alpha (\cos\alpha h^{3}_{11} -\cos\alpha h^{4}_{12})  \\
   & & - \partial_{1}\cos\alpha (J_{14,2} + \cos\alpha h^{3}_{12} -\cos\alpha h^{4}_{22})   \\
   &=& \partial_{2}\cos\alpha\cos\alpha(H^{3}-\frac{\partial_{2}\cos\alpha}{\sin\alpha}+\frac{J_{12,2}}{\sin\alpha}) \\
   & & - \partial_{1}\cos\alpha\cos\alpha(H^{4}-\frac{\partial_{1}\cos\alpha}{\sin\alpha}+\frac{J_{12,1}}{\sin\alpha}) \\
   & & - \partial_{1}\cos\alpha J_{14,2} \\
   &=& -\frac{\cos\alpha}{\sin{\alpha}}(1+\beta\frac{\sin^{2}{\alpha}}{\cos^{2}{\alpha}})(|\partial_{2}\cos\alpha|^{2}+|\partial_{1}\cos\alpha|^{2})  \\
   & & +\partial_{2}\cos\alpha\frac{J_{12,2}}{\tan\alpha}+\partial_{1}\cos\alpha\frac{J_{12,1}}{\tan\alpha} - \partial_{1}\cos\alpha J_{14,2} \\
   &=& -\frac{\cos^{2}\alpha+\beta\sin^{2}\alpha}{\cos\alpha\sin{\alpha}}( |\nabla\cos\alpha|^{2})  \\
   & & +\partial_{2}\cos\alpha\frac{J_{12,2}}{\tan\alpha}+\partial_{1}\cos\alpha(\frac{J_{12,1}}{\tan\alpha}-J_{14,2}).
\end{eqnarray*}
  
Thus, we have  
\begin{eqnarray*}
    \sin\alpha(H_{,1}^{4}+H_{,2}^{3}) &=& -\beta\frac{\sin^{2}\alpha}{\cos^{2}\alpha}\triangle\cos\alpha + \beta\frac{\sin^{2}\alpha(2+\beta\sin^{2}\alpha)}{\cos^{3}\alpha}|\nabla\alpha|^{2} \\
   & & -\beta\frac{\sin\alpha}{\cos^{2}\alpha}( \partial_{1}\cos\alpha(2\frac{J_{12,1}}{\tan\alpha}-J_{14,2}) \\
   & & +\partial_{2}\cos\alpha(2\frac{J_{12,2}}{\tan\alpha}+J_{13,2}) ). 
\end{eqnarray*}
  
By $(\ref{trieq})$, we obtain that 
\begin{eqnarray*}
  \triangle\cos\alpha &=& -\frac{\beta^{2}\sin^{4}\alpha+2\cos^{4}\alpha+2\beta\sin^{\alpha}\cos^{2}\alpha}{\cos^{3}\alpha}|\nabla\alpha|^{2} \\ 
   & & -\frac{2\cos\alpha}{\sin^{2}\alpha}((J_{12,1})^{2}+(J_{12,2})^{2}) \\
   & & +\frac{2(2\cos^{2}\alpha+\beta\sin^{2}\alpha)}{\cos\alpha\sin^{2}\alpha}(J_{12,1}\partial_{1}\cos\alpha + J_{12,2}\partial_{2}\cos\alpha)  \\
   & & -\beta\frac{\sin^{2}\alpha}{\cos^{2}\alpha}\triangle\cos\alpha + \beta\frac{\sin^{2}\alpha(2+\beta\sin^{2}\alpha)}{\cos^{3}\alpha}|\nabla\alpha|^{2} \\
   & & -\beta\frac{\sin\alpha}{\cos^{2}\alpha}( \partial_{1}\cos\alpha(2\frac{J_{12,1}}{\tan\alpha}-J_{14,2}) +\partial_{2}\cos\alpha(2\frac{J_{12,2}}{\tan\alpha}+J_{13,2}) ) \\
   & & -\sin\alpha(K_{1213}-K_{1224})  \\
   & & + J_{12,kk} + 2( -\beta\frac{\sin\alpha}{\cos^{2}\alpha}(J_{13,2}\partial_{2}\cos\alpha + J_{14,2}\partial_{1}\cos\alpha) ).
\end{eqnarray*}  
Therefore,
 \begin{eqnarray*}
    \triangle\cos\alpha &=& \frac{2\beta\sin^{2}\alpha}{\cos\alpha(\cos^{2}\alpha+\beta\sin^{2}\alpha)}|\nabla\alpha|^{2}-2\cos\alpha|\nabla\alpha|^{2}  \\
     & & -\frac{\sin\alpha\cos^{2}\alpha}{\cos^{2}\alpha+\beta\sin^{2}\alpha}(K_{1213}-K_{1224}) + \Theta   \\
     & & -\frac{\beta\sin\alpha}{\cos^{2}\alpha+\beta\sin^{2}\alpha}(\partial_{1}\cos\alpha J_{14,2}+3\partial_{2}\cos\alpha J_{13,2}), 
  \end{eqnarray*}
  where
  \begin{eqnarray*}
  \Theta &=& \frac{2\cos\alpha}{\sin^{2}\alpha}(1+\frac{\cos^{2}\alpha}{\cos^{2}\alpha+\beta\sin^{2}\alpha})(\partial_{1}\cos\alpha J_{12,1}+\partial_{2}\cos\alpha J_{12,2}) \\
      & & +\frac{\cos^{2}\alpha}{\cos^{2}\alpha+\beta\sin^{2}\alpha}J_{12,kk} - \frac{2\cos^{3}\alpha}{\sin^{2}\alpha(\cos^{2}\alpha+\beta\sin^{2}\alpha)}((J_{12,1})^{2}+(J_{12,2})^{2}).
  \end{eqnarray*}
 
\end{proof}

Observe the expression of $\Theta$, if we require $(\bar\nabla_{X}J(Y)+\bar\nabla_{Y}J(X))^{\top}=0$, we have $J^{k}_{i,i}=0$. Then $\Theta = 0$. Since $dim(\Sigma)=2$, this condition means that $(\nabla_{X} (J|_\Sigma))^{\top}=0$. Note that, if $\cos\alpha\equiv 1$, $J_{\Sigma}$ is the complex structure of $\Sigma$. So, the torsion tensor of $J$ is $N=0$, and $d\omega|_{\Sigma}=0$ means that $(\nabla_{X} (J|_\Sigma))^{\top} = \nabla_{X}(J|_\Sigma)=0$.


This paper gives the basic properties of $\beta$-symplectic critical surfaces. In addition, a natural problem is to construct a (almost-)Hermite manifold and the symplectic submanifold so that its complex structure $J$ satisfies the above conditions, i.e. $(\ref{Jc1})$ and $(\ref{Jc2})$.


\end{document}